\documentclass[11pt]{amsart}
\usepackage[all]{xy}
\usepackage{amssymb}
\usepackage{amsmath}
\usepackage{amsmath}

\swapnumbers
\newtheorem{thm}[subsection]{Theorem}
\newtheorem{prop}[subsection]{Proposition}
\newtheorem{lemma}[subsection]{Lemma}
\newtheorem{cor}[subsection]{Corollary}

\def\oo{\omega}
\def\dd{\delta}
\def\DD{\Delta}

\def\ss{\sigma}

\def\t{\otimes}

\def\Vect{{\texttt{Vect} }}

\def\pt{\textrm{-}}

\newcommand{\Vt}[1]{V^{\otimes #1}}

\def\row#1#2#3{(#1_{#2},\ldots,#1_{#3})}

\def\KK{{\mathbb{K}}}

\def\AA{{\mathcal{A}}}
\def\CC{{\mathcal{C}}}
\def\HH{{\mathcal{H}}}
\def\PP{{\mathcal{P}}}

\def\Prim{\mathrm{Prim\, }}
\def\Id{\mathrm{Id }}
\def\id{\mathrm{id }}

\def\Ker{\mathop{\rm Ker}}

\def\HHo{\overline {\HH}}

\def\Vtn{V^{\t n}}

\def\epi{\twoheadrightarrow}

\def\oleft#1{\omega^{#1}}
\def\oright#1{\overrightarrow{\omega^{#1}}}

\def\alg{\textrm{-alg}}
\def\bialg{\textrm{-bialg}}

\def\pt{\textrm{-}}

\def\row#1#2#3{(#1_{#2},\ldots,#1_{#3})}

\newenvironment{proo}{\begin{trivlist} \item{\emph{Proof.}}}
  {\hfill $\square$ \end{trivlist}}

\def\cpbtroisquatre{\vcenter{\xymatrix@R=6pt@C=6pt{
&&&&\\
&*{}\ar@{-}[u]&&*{}\ar@{-}[u]&\\
&&\mu\ar@{-}[ul]\ar@{-}[uu]\ar@{-}[ur]&&\\
&&\dd\ar@{-}[u]\ar@{-}[dll]\ar@{-}[dl] \ar@{-}[drr]\ar@{-}[dr] &&\\
*{}\ar@{-}[d]&*{}\ar@{-}[d]&&*{}\ar@{-}[d]&*{}\ar@{-}[d]\\
&&&&
}}}

\def\cpbmultiple{\vcenter{\xymatrix@R=6pt@C=6pt{
&&&&    &&&&    &&&&     \\
&&&&  *{}\ar@{-}[u]  &&&&   *{}\ar@{-}[u]   &&&&    \\
*{}\ar@{-}[rrrrrrrrrrr] &*{}\ar@{-}[uu]&&*{}\ar@{-}[ur]&  *{}\ar@{-}[uu]  &*{}\ar@{-}[ul]&*{}\ar@{-}[urr]  &*{}\ar@{-}[ur]  &  *{}\ar@{-}[u]    &&*{}\ar@{-}[ull]  &*{}&    \\
&*{}&*{}&*{}&  *{} &*{}& \omega &*{}& *{}&*{}&*{}&*{}& *{}   &*{}&*{}&*{}& *{}   \\
*{}\ar@{-}[uu] \ar@{-}[rrrrrrrrrrr] &*{}&*{}&*{}&   *{} &*{}&*{}&*{}&  *{}&*{}& *{}    &*{}\ar@{-}[uu] &  *{}  \\
&&*{}\ar@{-}[d]\ar@{-}[ul]\ar@{-}[ur]  &&  *{}\ar@{-}[d] \ar@{-}[u]
  &&*{}\ar@{-}[d] \ar@{-}[ul] \ar@{-}[u]\ar@{-}[ur]  &&&*{}\ar@{-}[d] \ar@{-}[ul] \ar@{-}[ur] &    \\
&&&&    &&&&      &&&&    \\
}}}

\def\cpbdeuxdeux{\vcenter{\xymatrix@R=2pt@C=2pt{
&&\\
*{}\ar@{-}[u]\ar@{-}[dr]&&*{}\ar@{-}[u]\ar@{-}[dl]\\
&*{}\ar@{-}[d]&\\
&*{}\ar@{-}[dl]\ar@{-}[dr]&\\
*{}\ar@{-}[d]&&*{}\ar@{-}[d]\\
&&
}}}

\def\cpbA{\vcenter{\xymatrix@R=2pt@C=2pt{
\ar@{-}[ddddd]&&\ar@{-}[ddddd]\\
&&\\
&&\\
&&\\
&&\\
&&
}}}

\def\cpbB{\vcenter{\xymatrix@R=2pt@C=2pt{
\ar@{-}[dd]&&\ar@{-}[dd]\\
&&\\
*{}\ar@{-}[drr]& &*{}\ar@{-}[dll]\\
*{}\ar@{-}[dd]& &*{}\ar@{-}[dd]\\
&&\\
&&
}}}

\def\cpbC{\vcenter{\xymatrix@R=2pt@C=2pt{
&&&&\\
&*{}\ar@{-}[dl]\ar@{-}[dr]\ar@{-}[u]&&*{}\ar@{-}[d]*{}\ar@{-}[u]&\\
*{}\ar@{-}[d]&&*{}\ar@{-}[dr]&*{}\ar@{-}[dl]&\\
*{}\ar@{-}[dr]&&*{}\ar@{-}[dl]&*{}\ar@{-}[dd]&\\
&*{}\ar@{-}[d]&&&\\
&&&&
}}}

\def\cpbD{\vcenter{\xymatrix@R=2pt@C=2pt{
&&&&\\
&*{}\ar@{-}[dl]\ar@{-}[dr]\ar@{-}[u]&&&\\
*{}\ar@{-}[ddd]&&*{}\ar@{-}[d]&&\\
&&*{}\ar@{-}[dr]&&*{}\ar@{-}[dl]*{}\ar@{-}[uuu]\\
&&&*{}\ar@{-}[d]&\\
&&&&
}}}

\def\cpbE{\vcenter{\xymatrix@R=2pt@C=2pt{
&&&&\\
&&&*{}\ar@{-}[u]\ar@{-}[dl]\ar@{-}[dr]&\\
&&*{}\ar@{-}[d]&&*{}\ar@{-}[ddd]\\
*{}\ar@{-}[uuu]\ar@{-}[dr]&&*{}\ar@{-}[dl]&&\\
&*{}\ar@{-}[d]&&&\\
&&&&
}}}

\def\cpbF{\vcenter{\xymatrix@R=2pt@C=2pt{
&&&&\\
&&*{}\ar@{-}[u]&&\\
*{}\ar@{-}[uu]*{}\ar@{-}[dr]&*{}\ar@{-}[dl]\ar@{-}[ur]&&*{}\ar@{-}[ul]\ar@{-}[d]\\
*{}\ar@{-}[dd]&*{}\ar@{-}[dr]&&*{}\ar@{-}[dl]\\
&&*{}\ar@{-}[d]&&\\
&&&&
}}}

\def\cpbG{\vcenter{\xymatrix@R=2pt@C=2pt{
&&&&&\\
&*{}\ar@{-}[dl]\ar@{-}[dr]\ar@{-}[u]&&&*{}\ar@{-}[dl]\ar@{-}[dr]\ar@{-}[u]&\\
*{}\ar@{-}[d]&&*{}\ar@{-}[dr]&*{}\ar@{-}[dl]&&*{}\ar@{-}[d]\\
*{}\ar@{-}[dr]&&*{}\ar@{-}[dl]&*{}\ar@{-}[dr]&&*{}\ar@{-}[dl]\\
&*{}\ar@{-}[d]&&&*{}\ar@{-}[d]&\\
&&&&&&
}}}

\def\arbreA{\vcenter{\xymatrix@R=3pt@C=3pt{
&& \\
&*{}\ar@{-}[ur] \ar@{-}[ul] \ar@{-}[d]     &\\
&&
}}}

\def\arbreAgrand{\vcenter{\xymatrix@R=30pt@C=30pt{
&& \\
&*{}\ar@{-}[ur] \ar@{-}[ul] \ar@{-}[d]     &\\
&&
}}}

\def\arbreBA{\vcenter{\xymatrix@R=2pt@C=2pt{
&&&&\\
&&&*{}\ar@{-}[ul] & \\
&&*{}\ar@{-}[uurr] \ar@{-}[uull] \ar@{-}[d]     &&\\
&&&&
}}}

\def\arbreAB{\vcenter{\xymatrix@R=2pt@C=2pt{
&&&&\\
&*{}\ar@{-}[ur] &&& \\
&&*{}\ar@{-}[uurr] \ar@{-}[uull] \ar@{-}[d]     &&\\
&&&&
}}}

\def\arbreABDA{\vcenter{\xymatrix@R=1pt@C=1pt{
1&&2&&3&&4&&5\\
&&&& &&&&\\
&*{}\ar@{-}[ur] &&&&&&*{}\ar@{-}[ul] & \\
&&*{}\ar@{-}[uurr] &&&&&&\\
&&&&&&&&\\
&&&&*{}\ar@{-}[uuuurrrr] \ar@{-}[uuuullll] \ar@{-}[d] &&&&\\
&&&& &&&&
}}}

\begin{document}

\author[R. Holtkamp, J.-L. Loday, M. Ronco]{Ralf Holtkamp, Jean-Louis Loday, Mar\'ia Ronco}
\address{RH: Fakult\"at f\"ur Mathematik, Ruhr-Universit\"at,
          44780 Bochum, Germany}
\email{ralf.holtkamp@ruhr-uni-bochum.de}
\address{JLL: Institut de Recherche Math\'ematique Avanc\'ee,
    CNRS et Universit\'e Louis Pasteur,
    7 rue R. Descartes,
   67084 Strasbourg Cedex, France}
\email{loday@math.u-strasbg.fr}
\address{MOR: Departamento de Matematicas,
  Facultad de Ciencias-Universidad de Valparaiso,
      Avda. Gran Bretana 1091,
         Valparaiso, Chile}
\email{maria.ronco@uv.cl}

\title{Coassociative magmatic bialgebras and the Fine numbers}
\subjclass[2000]{16A24, 16W30, 17A30, 18D50, 81R60.}
\keywords{ Bialgebra, Hopf algebra, Cartier-Milnor-Moore,
Poincar\'e-Birkhoff-Witt, magmatic, operad, Fine number, pre-Lie algebra}

\date{\today}

\begin{abstract}
We prove a structure theorem for the connected coassociative magmatic bialgebras. The space of primitive elements is an algebra over an operad called the primitive operad. 
We prove that the primitive operad is magmatic generated by $n-2$ operations of arity $n$. The dimension of the space of all the $n$-ary operations of this primitive operad 
turns out to be the Fine number $F_{n-1}$. In short, the triple of operads $(As, Mag, MagFine)$ is good.

\end{abstract} \maketitle

\section*{Introduction} \label{S:int} A magmatic algebra is a vector space equipped with a unital binary operation, denoted $x\cdot y$  with no further assumption. Let us suppose that it is also equipped with a counital binary cooperation, denoted $\DD(x)$. There are different compatibility relations that one can suppose between the operation and the cooperation. Here are three of them:

$$\begin{array}{lrcl}
\textrm{Hopf:}                      &\DD(x\cdot y) & = & \DD(x)\cdot \DD(y)\\
\textrm{magmatic:}               &  \DD(x\cdot y) & = &x\cdot y\t 1 +  x\t y+ 1\t x\cdot y\\
\textrm{unital infinitesimal:}  &\DD(x\cdot y) & = & \DD(x) \cdot (1\t y) + (x\t 1)\cdot \DD(y) - x\t y\\
\end{array}$$

On the free magmatic algebra the cooperation $\DD$ is completely determined by this choice of compatibility relation and the assumption that the generators are primitive.
 In the Hopf case $\DD$ is coassociative and cocommutative, giving rise to the notion of $Com^c\textrm{-}Mag$-bialgebra. This case has been
addressed in \cite{H2}.
 In the magmatic case $\DD$ is comagmatic, giving rise to the notion of $Mag^c\textrm{-}Mag$-bialgebra. This
case has been addressed by E.~Burgunder in \cite{B}.  In the u.i.~case, $\DD$ is coassociative, giving rise to the notion of
$As^c\textrm{-}Mag$-bialgebra. This case is the subject of this paper.

First, we construct a functor $F:Mag\alg \to MagFine\alg$, where $MagFine$ is the operad of algebras 
having $(n-2)$ operations of arity $n$ for $n\geq 2$ (no binary operation, one ternary operation, etc.). Then we show that the primitive part of any 
$As^c\textrm{-}Mag$-bialgebra is closed under the $MagFine$ operations. The functor $F$ has a left adjoint denoted $U:MagFine\alg \to Mag\alg$.

Second, we prove a structure theorem for connected $As^c\textrm{-}Mag$-bialgebras. It says that, as a bialgebra, any connected $As^c\textrm{-}Mag$-bialgebra $\HH$ is 
isomorphic to $U(\Prim \HH)$ and that, as a coalgebra, it is cofree. Observe that $U(\Prim \HH)$ has a meaning since the primitive part is a $MagFine$-algebra. 
These statements are analogous to the Cartier-Milnor-Moore theorem (cf. \cite{Ca, MM}) and the Poincar\'e-Birkhoff-Witt theorem for cocommutative (classical) 
bialgebras.

As a consequence there is an equivalence of categories:
\begin{displaymath} \{\textrm{connected\ } As^c\pt Mag\bialg \}
{\begin{array}{c} U\\
\leftrightarrows\\
\Prim
\end{array}}
\{MagFine\alg \} \ .
\end{displaymath}

 In the terminology of \cite{Lod}, the triple of operads $$(As, Mag, MagFine)$$ is a good triple of operads.

As a byproduct of some intermediate Lemma we obtain a new combinatorial interpretation of the Fine numbers $F_{n-1},
n=1,2,\ldots$ (cf.~\ref{cor:Fine}): $$(1, 0, 1, 2, 6, 18, 57,186, 622, 2120, 7338, 25724, \ldots   ).$$

In the magmatic case, $(Mag, Mag, Vect)$ is a good triple of operads, see \cite{B}.

The study of primitive elements in the Hopf case, cf.~\cite{HS, SU}, leads to a good triple $(Com, Mag, Sabinin)$, due to
Shestakov and others. Here the integer series $\frac{1}{(n-1)!}{\dim Sabinin(n)}$ is given by the Log-Catalan numbers
$$(1, 1, 4, 13, 46, 166, 610, 2269, 8518, 32206, \ldots   ).$$

For other instances of strong relationship between bialgebras and trees (generators of the free magmatic algebra), the reader can look at (nonexhaustive list) \cite{A, B, H1,
H2, LR0, LR1, LR2, R}.

 \bigskip

\noindent {\bf Notation.} In this paper $\KK$ is a field and all vector spaces are over $\KK$. Its unit is denoted $1_\KK$ or just 1. 
The category of vector spaces over $\KK$ is denoted by $\Vect$. The vector space spanned by the elements of a set $X$ is denoted $\KK[X]$. 
The tensor product of vector spaces over $\KK$ is denoted by $\t$. The tensor product of $n$ copies of the space $V$ is
denoted by $\Vt n$. For $v_i\in V$ the element $v_1\t \cdots \t v_n$ of $\Vt n$ is denoted by  $(v_1, \ldots , v_n)$ or
simply by $v_1 \ldots v_n$.
A linear map $\Vt n \to V$ is called an \emph{ $n$-ary operation} on $V$ and a linear map $ V\to \Vt n$ is called an \emph{
$n$-ary cooperation} on $V$.  The symmetric group is the automorphism group of the finite set $\{1,\ldots , n\}$ and is denoted $S_n$. 
It acts on $\Vtn$ by $\ss \cdot \row v1n = \row v{\ss^{-1}(1)}{\ss^{-1}(n)}$. The action is extended to an action of $\KK[S_n]$ by linearity.

For a given type of algebras, or algebraic operad $\PP$,
 the category of algebras of type $\PP$, also called $\PP$-algebras, is denoted $\PP\alg$. All the operads $\PP$ appearing in this paper are \emph{regular operads}. 
So, the space $\PP(n)$ of all $n$-ary operations is of the form $\PP(n)= \PP_n \t \KK[S_n]$. The free algebra of type $\PP$ over the vector space $V$  is $\PP(V)= \bigoplus_{n\geq 0}\PP_n\t V^{\t n}$. We often look at $\PP$ as an endofunctor of the category $\Vect$ of vector spaces. For more information on algebraic operads see for instance the first section of \cite {Lod}.

When $\PP=As$, the operad of unital associative algebras, the free algebra over $V$ is the tensor algebra $T(V)$. We use equivalently the notation $As(V)$ and $T(V)$.

\section{Coassociative magmatic bialgebras}

In this section we introduce the notion of coassociative magmatic bialgebra and we prove that a free magmatic algebra has a natural  structure of $As^c\textrm{-}Mag$-bialgebra. We give an explicit description of the coproduct.

\subsection{Definition}\label{magbialg} A \emph{coassociative magmatic bialgebra} $\HH$, or $As^c\textrm{-}Mag$-bialgebra for short, is a vector space $\HH$ over $\KK$ equipped with a coassociative cooperation $\DD: \HH\to \HH\t \HH, x\mapsto \sum x_{(1)}\t x_{(2)},$ and  a binary operation $(x,y)\mapsto x\cdot y$, which satisfy the \emph{unital infinitesimal relation}:
\begin{eqnarray*}
\DD(x\cdot y) &=& \DD(x) \cdot (1\t y) + (x\t 1)\cdot \DD(y) - x\t y\ , \ i.e.\ \\
\sum (x\cdot y)_{(1)}\t (x\cdot y)_{(2)} &=& \sum x_{(1)}\t x_{(2)}\cdot y + \sum x\cdot y_{(1)}\t y_{(2)} - x\t y\ .\\
\end{eqnarray*}
Moreover we suppose that the operation $\cdot$ has a unit $1$, that the cooperation $\DD$ has a counit and that $\DD(1)=1\t 1$. We denote by $\HHo$ the augmentation ideal: $\HH = \KK\ 1 \oplus \HHo$.

Observe that we do not suppose anything else about the operation, in particular it is not supposed to be associative, hence it is a \emph{magmatic product}.

It will be helpful to work with the \emph{reduced coproduct} $\dd : \HHo \to \HHo\t \HHo$ which is defined as
$$ \dd(x) = \DD(x) - x\t 1 - 1\t x.$$

The \emph{compatibility relation} satisfied by the reduced coproduct is
$$\dd(x\cdot y) = \dd(x) \cdot (1\t y) + (x\t 1)\cdot \dd(y) + x\t y\ ,$$
or equivalently
$$\dd(x\cdot y) = \sum x_{(1)}\t x_{(2)}\cdot y + \sum  x \cdot y_{(1)}\t y_{2} + x\t y\ ,$$
where $\dd (x)=\sum x_{(1)}\otimes x_{(2)}$.
It is called the \emph{unital infinitesimal} compatibility relation (for the nonunital framework) and it can be pictured as:
$$\cpbdeuxdeux \ =\   \cpbD\ +\ \cpbE\ +\ \cpbA  \ .$$

We say that an $As^c\textrm{-}Mag$-bialgebra $\HH$ is {\it connected}, or conilpotent,  if $\HH = \bigcup _{r\geq 0} F_r\HH$
where $F_r\HH$ is the coradical filtration of $\HH$ defined recursively by the formulas
\begin{eqnarray*}
F_0\HH &:=& \KK\ 1, \\
F_r\HH & :=& \{ x\in \HH \mid  \dd  (x) \in F_{r-1}\HH \t F_{r-1}\HH \}\ .
\end{eqnarray*}

By definition the \emph{space of primitive elements} is defined as $$\Prim \HH := \Ker \dd \subset \HHo.$$

\subsection{Free magmatic algebra} Let us recall that the free magmatic algebra over $V$ is of the form
$$Mag(V) =\bigoplus_{n\geq 0} \KK[Y_{n-1}]\t \Vtn$$
where $Y_{n-1}$ is the set of planar binary rooted trees with $n$ leaves (and $ \KK[Y_{-1}]\t V^{\t 0}= \KK\ 1$) . The magmatic product is induced by the grafting of trees and the concatenation of tensors:
$$(t;v_1\ldots v_p)\cdot (s;v_{p+1}\ldots v_{p+q}) = (t\vee s;v_{1}\ldots v_{p+q}),$$
see for instance \cite{H1}.

\begin{prop}\label{MagV} The free magmatic algebra $Mag(V)$ over the vector space $V$ has a natural structure of coassociative magmatic bialgebra.
\end{prop}

\begin{proo}
Let us equip the tensor product $Mag(V)\t Mag(V)$ with the following magmatic product:
$$(x\t y)\cdot (x'\t y') := (x\cdot x')\t (y\cdot y').$$
Since $Mag(V)$ is free, the map $V\to Mag(V)\t Mag(V), v\mapsto v\t 1 + 1\t v$ admits a unique lifting $\DD:Mag(V)\to Mag(V)\t Mag(V)$ which satisfies the unital infinitesimal relation. Observe that $\DD$ is not a magmatic homomorphism.

The map $\DD$ is coassociative because the two maps $(\DD\t\Id)\DD$ and $(\Id \t \DD) \DD$ extend the same map
 $Mag(V)\to Mag(V)^{\t 3}, v\mapsto v\t 1\t 1 + 1\t v\t 1 + 1 \t 1 \t v$ under the compatibility relation. Hence they are equal on $Mag(V)$.
\end{proo}

We denote by $\Prim Mag(V)$ the primitive part of $Mag(V)$ under this cooperation.

\begin{cor} The functor $V\mapsto \Prim Mag(V)$ is an algebraic operad under the composition of the operad $Mag$. In other words there is a natural transformation of functors 
$\Prim Mag \circ \Prim Mag\ \cdots\!\!\! >  \Prim Mag$, which makes the following diagram commutative
$$\xymatrix{
\Prim Mag \circ \Prim Mag  \ar@{.>}[r] \ar@{>->}[d] & \Prim Mag\ar@{>->}[d]\\
Mag \circ  Mag          \ar[r]                                              & Mag
}$$ 

The primitive part of a coassociative magmatic
bialgebra is  a $\Prim Mag$-algebra.
\end{cor}

\begin{proo} An $As^c\textrm{-}Mag$-bialgebra is a generalized bialgebra in the sense of
\cite{Lod} because the compatibility relation between the operation and the cooperation is distributive.

In \cite{Lod} it is proved that, for a generalized bialgebra of type $\CC^c\pt \AA$, if the free algebra $\AA(V)$ is a $\CC^c\pt \AA$-bialgebra, then 
$\Prim \AA(V)$ determines an operad. By Proposition \ref{MagV} this hypothesis is fulfilled in the $As^c\textrm{-}Mag$-case, whence the assertion.
\end{proo}

\subsection{The coproduct $\dd$ made explicit}\label{Properties}
 Let $\HH$ be an $As^c\textrm{-}Mag$-bialgebra over the field $\KK$. Let $\oo^n$ be the operation which corresponds to the left comb, that is
$$\oo^n(x_{1},\ldots ,x_{n}):= ((\ldots ((x_1\cdot x_2)\cdot x_3)\ldots)\cdot x_n ) .$$
Given primitive elements $x_1,\ldots,x_n$ of $\HH$, it is immediate to check that:
$\dd(\oleft{n}(x_1,\ldots,x_n)) = x_1\t
\oleft{n-1}(x_2,\ldots,x_n) +\cdots\hfill $\\
${}\qquad + \oleft{i}(x_1,\ldots,x_{i})\t  \oleft{n-i}(x_{i+1},\ldots,x_{n}) +\cdots+\oleft{n-1}(x_1,\ldots,x_{n-1})\t x_n\ .
$

By Proposition \ref{MagV}, the free magmatic algebra is a coassociative magmatic bialgebra, where the product is induced by the grafting of trees.

We now give an explicit description of $\dd$ in terms of split trees. Let $t$ be a planar binary rooted tree with $r$ leaves numbered from left to right by 
$1, 2, \ldots , r$. Let $i,\ 1\leq i< r$ be an integer. We split the tree $t$ into two trees
$t_{(1)}^{i}$ and $t_{(2)}^{i}$ by cutting in between the leaves
$i$ and $i+1$. More exactly, the tree $t_{(1)}$ is the part of $t$
 which is on the left side of the path from leaf $i$ to the root
(including the path). The tree $t_{(2)}$ is the analogous part on
the right side of the path which runs from leaf $i+1$ to the
root. For instance,
$$ t= \arbreABDA ,\  i=2 ,\  t_{(1)}^{2}= \arbreA ,\   t_{(2)}^{2}= \arbreBA\ .$$

\begin{lemma}\label{splittinglemma}  Let $Mag(V)$ be the free magmatic algebra on $V$.
Then the linear map $\dd$ applied to an element $(t;v_{1}\ldots
v_{r})$ can be written as a sum $\sum_{i=1}^{r-1} \dd_i$, where
$$\dd _{i} (t ;v_{1}\ldots v_{r})= (t_{(1)}^{i};v_{1}\ldots v_{i})\t
(t_{(2)}^{i};v_{i+1}\ldots v_{r}).$$
\end{lemma}

\begin{proo} The result is immediate for $r=1,2$. For $r>2$, there exist unique trees $t^{l}$ and $t^{r}$ such that $t= t^{l}\vee t^{r}$,
with $t^{l}\in Y_{k-1}$ and $t^{r}\in Y_{r-k-1}$ for some $1\leq
k\leq r-1$. Applying the recursive hypothesis and the relation between the product
$\cdot $ and the coproduct $\dd$, one gets: $$\displaylines { \dd
(t;v_{1}\ldots v_{r})=\dd ((t^{l};v_{1}\ldots v_{k})\cdot
(t^{r};v_{k+1}\ldots v_{r}))=\cr
 \sum_{i=1}^{k-1}\big((t^{l})_{(1)}^{i};v_{1}\ldots v_{i}\big)\t
\big((t^{l})_{(2)}^{i}\vee t^{r};v_{i+1}\ldots v_{r})+ (t^{l};v_{1}
\ldots v_{k})\t (t^{r};v_{k+1}\ldots v_{r}\big)+\cr
 \sum_{i=k+1}^{r-1}\big(t^{l}\vee (t^{r})_{(1)}^{i-k};v_{1}\ldots v_{i}\big)\t
\big((t^{r})_{(2)}^{i-k};v_{i+1}\ldots v_{r}\big).\cr }$$

This formula implies that $\dd$ satisfies the u.i.  compatibility relation, therefore it is the (reduced) magmatic coproduct of Proposition \ref{MagV}.
 \end{proo}

\section{MagFine algebras}

We introduce the notion of $MagFine$-algebra and we compute $\dim
MagFine_n$. We mention a new combinatorial interpretation of the
Fine numbers. We construct a functor $Mag\alg \to MagFine\alg.$ We
will later show that it factors surjectively through $\Prim
Mag\alg$.

\subsection{$MagFine$ algebras} By definition a \emph{$MagFine$ algebra} is a vector space endowed with $(n-2)\  n$-ary operations denoted $m^n_i$, $1\leq i\leq n-2$,  for all
$n\geq 3$. Let $MagFine$ be the operad associated to this type of algebras. Since the generating operations are not supposed to satisfy any symmetry
property nor any relations, the operad is regular. In particular the space of $n$-ary operations is of the form
$$MagFine(n) = MagFine_n\t \KK[S_n]\ .$$

\begin{prop}\label{Fineseries} The dimension of $MagFine_n$ is the Fine number $F_{n-1}$ defined by the generating series
$$F(x)=\sum_{n\geq 1} F_{n-1} x^n = \frac{1+2x-\sqrt{1-4x}}{2(2+x)}\
.$$
\end{prop}

\begin{proo} Let $\PP$ be the operad defined by $(n-2)$ $n$-ary operations and whose relations are: any nontrivial composition is zero. This is a nilpotent operad, so its Koszul dual is free with the same generators. Hence $\PP^! = MagFine$.

We will use a result of B. Vallette \cite{V} which relates the dimensions of $\PP(n)$ and of $\PP^!(n)$.

For a presented quadratic  operad $\PP$ we let $\PP_{(d)}(n)$ be the space spanned by the $n$-ary operations constructed out of  $d$ generating operations. 
In the case at hand we have $\dim \PP_{(0)}(1) =1, \dim \PP_{(1)}(n) =(n-2) n! $ and $\dim \PP_{(d)}(n) =0$ otherwise.

It is proved in \cite{V} Section 9, that the Poincar\'e series
$$f_{\PP}(x,z):= \sum_{{n\ge 1}\atop {d\geq 0}}\frac {\dim  \PP_{(d)}(n)}{n!} x^n z^d $$
 of $\PP$ and $f_{\PP^!}(x,z)$ are related by the relation
$$f_{\PP^!}(f_{\PP}(x,z),-z)=x\ ,$$
when $\PP$ is a Koszul operad. Since $\PP^{!}$ is a free operad, it is Koszul and we can apply this formula.
We know that, in our case, 
$$f_{\PP}(x,z)= x + z\sum_{n\geq
3}(n-2)x^n =  x + z\frac{x^{3}}{(1-x)^2}\ .$$
 We want to compute the
numbers $\dim MagFine_{n}$, or equivalently the series 
$$F(x)=
\sum_{n\geq 1}\dim MagFine_{n}x^n =\sum_{n\geq 1}\frac{ \dim
MagFine(n)}{n!}x^n  = f_{\PP^!}(x,1).$$
 Applying Vallette's formula
for $z=-1$ we get
$$F\big(x - \frac{x^3}{(1-x)^2}\big) = x .$$
It follows that $F(x) = \frac{1+2x - \sqrt{1-4x}}{2(2+x)} $. This is exactly the generating function of the Fine numbers, cf.~\cite{DS}.
\end{proo}

\subsection{Remark} The Fine numbers  are, in low dimension (starting at $F_0=1$),

$$ 1, 0, 1, 2, 6, 18, 57,186, 622, 2120, 7338, 25724, \ldots  .$$

They admit several combinatorial interpretations,  cf.~\cite{DS}. Here is a new one, which is a consequence of Proposition \ref{Fineseries}.
\begin{cor}\label{cor:Fine} Let $X_{n}$ be the set of planar rooted trees with $n$ leaves, with no unary nor binary vertices and whose $k$-ary vertices are labelled by $\{1, \ldots, k-2\}$. Then $\# X_{n}=F_{n-1}$ (Fine number).
\end{cor}
\begin{proo} It is immediate to check that the trees in $X_{n}$ are in bijection with a basis of $MagFine_{n}$.
\end{proo}

\subsection{Higher operations on a magmatic algebra} The \emph{associator} operation $as$  is a ternary operation defined as
$$ as(x_1,x_2,x_3):=(x_1\cdot x_2)\cdot x_3-x_1\cdot (x_2\cdot x_3).$$
If $(A,\cdot)$ is any magmatic algebra over $\KK$, we can equip it with the $n$-ary operations $\mu_i^n: A^{\t n}\to A$,
for all $n\geq 3$, all $1\leq i \leq n-2$, defined as follows:
\begin{eqnarray*}
\mu_1^3(x_1,\ldots,x_3) &:=& as(x_1,x_2,x_3),
\\
\mu_1^n(x_1,\ldots,x_n) &:=&
as(x_1,\oleft{n-2}(x_2,\ldots,x_{n-1}),x_n),
\\
\mu_2^4(x_1,\ldots,x_4) &:=& as(x_1,x_2,x_3\cdot
x_4)-as(x_1,x_2,x_3)\cdot x_4\\ \mu_i^n(x_1,\ldots,x_n) &:=&
\mu_2^4(x_1,\oleft{n-i-1}(x_2,\ldots,x_{n-i}),\oleft{i-1}(x_{n-i+1},\ldots,x_{n-1}),x_n),
\end{eqnarray*}
for $x_1,\ldots,x_n\in A$, $n\geq 4$, and $2\leq i \leq n-2$. We recall that $\oo^n $ is the left comb product.

In other words we have constructed a morphism of operads
$$MagFine \to Mag, \quad m_i^n\mapsto \mu_i^n.$$
In the following, we show that it factors through $\Prim Mag$.
\medskip

\subsection{Associator} The associator operation $as$
is a $\Prim Mag$-operation: $as\in \Prim Mag(3)$. In other words, for any $As^c\textrm{-}Mag$-bialgebra
$\HH$, and primitive elements $x_1,x_2,x_3$, the element
$as(x_1,x_2,x_3)$ is again primitive. This is a special case of
the following Lemma:

\begin{lemma}\label{coprodAs} Let $x \in \Prim\HH$, $y,z\in \HH$. Then:
$$ \dd (as(x,y,z))=\sum as(x,y,z_{(1)})\t\ z_{(2)}.$$
\end{lemma}
\begin{proo}
We use the compatibility relation given in
\ref{magbialg}. Since $x \in \Prim\HH$, $\dd(x\cdot y)=(x\t
1)\cdot \dd(y) + x\t y=\sum x\cdot y_{(1)}\t y_{(2)} + x\t y$, and
$\dd (as(x,y,z))$ is given by
\begin{eqnarray*}
 &\sum x\cdot y_{(1)}\t y_{(2)}\cdot z + x\t
y\cdot z + \sum (x\cdot y)\cdot z_{(1)}\t z_{(2)} +x\cdot y\t z \\
 &-\sum
x\cdot y_{(1)}\t y_{(2)}\cdot z -\sum x\cdot (y\cdot z_{(1)})\t
z_{(2)}
 -x\cdot y\t  z
-x\t y\cdot z\\
 & = \sum as(x,y,z_{(1)})\t\ z_{(2)}. \\
\end{eqnarray*}
\end{proo}

\begin{prop}\label{MuOperations}
If $\HH$ is an $As^c\textrm{-}Mag$-bialgebra, then $\Prim\HH$ is closed under the operations $\mu_i^n$, for $n\geq 3$ and $1\leq i\leq n-2$ (in other words the operations $\mu_i^n$ lie in $\Prim Mag(n)$).
\end{prop}

\begin{proo}
For all $n\geq 3$, $2\leq i \leq n-2$, and $x_1,\ldots,x_n\in \Prim\HH$, we have to show that $\mu_i^n(x_1,\ldots,x_n)\in
\Prim\HH$, too.

Using Lemma \ref{coprodAs} in the case where $z=x_n\in \Prim\HH$, we immediately get that $\mu_1^n(x_1,\ldots,x_n)\in \Prim\HH$ for
all $n$.

We need to compute $\dd(\mu_i^n((x_1,\ldots,x_n))$, which is, in short notation,
$$\dd(as(x_1,\oleft{n-i-1},\oleft{i}))-\dd(as(x_1,\oleft{n-i-1},\oleft{i-1})\cdot
x_n).$$

 By Lemma \ref{coprodAs} (together with \ref{Properties}) it follows that
$$\dd(as(x_1,\oleft{n-i-1},\oleft{i}))=\sum_{j=1}^{i-1}
as(x_1,\oleft{n-i-1},\oleft{i-j})\t \oleft{j}.$$

Similarly,
\begin{eqnarray*}
&&\dd(as(x_1,\oleft{n-i-1},\oleft{i-1})\cdot x_n)\\
&=&\sum_{j=1}^{i-2} as(x_1,\oleft{n-i-1},\oleft{i-j})\t
\oleft{j}\cdot x_n + as(x_1,\oleft{n-i-1},\oleft{i-1})\t x_n
\\
&=&\sum_{j=1}^{i-1} as(x_1,\oleft{n-i-1},\oleft{i-j})\t \oleft{j},
\end{eqnarray*}
which ends the proof.
\end{proo}

\begin{prop}\label{injective} The morphism of operads
$$MagFine \to \Prim Mag,\qquad  m_i^n\mapsto \mu_i^n$$
 is injective.
\end{prop}
For $V$ a vector space, let $PMag(V)$ be the smallest subspace of
$Mag(V)$ which contains $V$ and is closed under all operations
$\mu_i^n$. All elements of $PMag(V)$ are primitive, by Proposition
\ref{MuOperations}. They can be presented by (linear combinations
of) compositions of operations $\mu_i^n$, evaluated on elements of
$V$, and we are going to show that this presentation is unique.

Before Proposition \ref{injective} we prove the following result.

\begin{lemma}\label{IntersectPMag}
For any $r\geq 1$, let $PMag^r$ be the subspace of $Mag(V)$ spanned
by elements of type $\oleft{r}(z_1,\ldots,z_r)$, with $z_i\in
PMag(V),$ and let similarly $PMag^{\neq r}= + _{{n\geq 1}\atop {n\neq r}}PMag^{n}$ be the space spanned by
all elements $\oleft{j}(z_1,\ldots,z_j)$, $j\neq r$. 
Then \begin{enumerate}\item $PMag^r\cap
PMag^{\neq r}=0$, for $r\geq 1$.
\item For any pair of positive integers $r,s$, the intersection
$$ (\bigoplus_{(i,j)\neq (r,s)} PMag^i\cdot PMag^j) \bigcap (PMag^r\cdot PMag^s)$$
 is zero.
\item The intersection
 $$(\bigoplus_{j \geq 1} PMag^j)\bigcap
 (\bigoplus_{r\geq 1, s\geq 2} PMag^r\cdot PMag^s)$$
is zero.
\end{enumerate}
\end{lemma}

\begin{proo}

\begin{enumerate} \item By the properties (\ref{Properties}) of $\dd$, we have that
$$\dd^{l-1}(\oleft{j}(z_1,\ldots,z_j))=\begin{cases} z_1\t\ldots\t
z_j, & \text{for } l=j\\ 0, & \text{for } l>j.\end{cases}$$

Since $Mag(V)$ is free, the restriction $\dd^{l-1}\vert_{PMag^l}:
PMag^l\to PMag^{\t l}$ is bijective.

 Let $\sum_{j=1}^n x_j=0$, with $x_j\in PMag^j.$ Applying $\dd^{n-1}$, we get that
$\sum_{j=1}^{n-1} x_j=x_n=0.$ The result follows applying a recursive argument on $n-i$.

\item Note that the first point of the Lemma implies that, for any subspace $S$ of $Mag(V)$ and any $r\geq 2$,
\begin{eqnarray*}
 (\bigoplus_{j=1}^{r-1}S\cdot PMag^j)\bigcap
(S\cdot PMag^r) &= & 0,\\
 (\bigoplus_{j=1}^{r-1}PMag^j\cdot S)\bigcap
(PMag^r\cdot S) &= & 0.\\
\end{eqnarray*}

Suppose that, for some $x\in PMag^r$ and some $y\in PMag^s$, there exist $\{ w_k\cdot z_k \}_{1\leq k \leq p}$ such
that:
 $$ x\cdot y = \sum_{k=1}^p w_k \cdot z_k,$$ 
 with $w_k\cdot z_k \in \bigoplus_{(i,j)\neq (r,s)} PMag^i\cdot PMag^j$. Since
$Mag(V)$ is free, there exists some $\{ z_{k_l} \}$ such that $y=\sum_l z_{k_l}$. The remark above implies that
 $z_{k_l} \in PMag^s$ for all $j$.

Then $x\cdot y=\sum_l w_{k_l}\cdot z_{k_l}$, which implies that $x=\sum_l w_{k_l}$, and, by the first point, $\sum_l
w_{k_l}\in PMag^r$. This proves the assertion.

\item The proof is immediate, because $Mag(V)$ is free, and the unique way to write an element of $z\in PMag^j$, $j\geq 1$, as
a product $z=x\cdot y$ is for $(x,y)=(1_{\KK},z)$ or $(x,y)=(z,1_{\KK}).$
 (We note that the case $s=1$ is excluded, because $PMag^s\cdot PMag^1=PMag^{s+1}.)$
 \end{enumerate}
\end{proo}

\medskip

\begin{proo} (Proposition \ref{injective})

Let $\mu_k^n(PMag^{\t n})$ be the image of the restriction of $\mu_k^n$ to $PMag(V)^{\t n}$, for $n\geq 3$, $1\leq k\leq n-2$.

\noindent We are going to prove that, for any $n$ and any $1\leq k\leq n-2$, the intersection 
$$ (\bigoplus_{(m,j)\neq (n,k)} \mu_j^m(PMag^{\t m})) \bigcap \mu_k^n(PMag^{\t n}) $$
 is zero. This implies that $PMag(V)=V\oplus (\bigoplus_{n,k}\ \mu_k^n(PMag^{\t n})).$

Let us look at the elements of $Mag(V)\cdot PMag^{1}$ which appear as terms in $\mu_k^n(PMag^{\t n}).$ We have to consider the
following cases:
\begin{itemize}
\item[(1)]
 In $\mu_1^n(x_1,\ldots,x_n)$ there is only one term of this type, the element $(x_1\cdot \oleft{n-2}(x_2,\ldots,x_{n-1}))\cdot
x_n$, for $x_1,\ldots,x_n\in PMag^{1}.$
\item[(2)]
In $\mu_k^n(x_1,\ldots,x_n)$, for $k\geq 2$, there are two terms of this type, the element $((x_1\cdot
\oleft{n-k-1}(x_2,\ldots,x_{n-k}))\cdot \oleft{k-1}(x_{n-k+1},\ldots,x_{n-1}))\cdot x_n$ and $(x_1\cdot
(\oleft{n-k-1}(x_2,\ldots,x_{n-k})\cdot \oleft{k-1}(x_{n-k+1},\ldots,x_{n-1})))\cdot x_n,$ for
$x_1,\ldots,x_n\in PMag^{1}.$
\end{itemize}

While an element of type $x_1\cdot\oleft{n-2}(x_2,\ldots,x_{n-1})$ is in
$$PMag^{1}\cdot PMag^{n-2}, \text{ for } n\geq 3,$$
 an element of type $x_1\cdot (\oleft{n-k-1}(x_2,\ldots,x_{n-k})\cdot
\oleft{k-1}(x_{n-k+1},\ldots,x_{n-1}))$ belongs to $$PMag^{1}\cdot (PMag^{n-k-1}\cdot PMag^{k-1}).$$

We thus note that:
\begin{itemize}
\item[(1)]
For $r\geq 1$ and $s\geq 2$, the unique terms of
$$(PMag^{1}\cdot
(PMag^r\cdot PMag^s))\cdot PMag^{1}
$$
which appear in $PMag(V)$ are in the image of $\mu_{s+1}^{r+s+2}.$

\item[(2)]
For $k\geq 1,$ the unique terms of $$ (PMag^{1}\cdot PMag^{k})\cdot 
PMag^{1}$$ which appear in $PMag(V)$ are in the image of 
$\mu_1^{k+2}$ or $\mu_2^{k+2}.$ 
In this case we note that terms of $((PMag^1\cdot PMag^{k-1})\cdot 
PMag^1)\cdot PMag^1$, for $k\geq 2$, characterize the image of 
$\mu_2^{k+2}$, whereas $$((PMag^1\cdot PMag^{k-1})\cdot 
PMag^1)\cdot PMag^1\cap\mu_1^{k+2}(PMag^{\t k+2})=0.$$
\end{itemize}

We get the result by an application of Lemma \ref{IntersectPMag}.

The results above prove that $Prim(Mag(V))$ is freely generated by the operations $\mu _{i}^{n}$, with $n\geq 3$ and $1\leq i\leq n-2$. 
So, we may conclude that the operads $\Prim Mag $ and $MagFine $ coincide.

\end{proo}

\section{On the structure of  coassociative magmatic bialgebras}

 We constructed a functor $$F:Mag\alg \to MagFine\alg ,\  F(A, \cdot) = (A,\mu^{n}_{i}). $$
We denote by $U:MagFine\alg \to Mag\alg$ the left adjoint
functor of $F$. So, if A is a $MagFine$-algebra, then $U(A)$ is the quotient of $Mag(A)$ by the relations which consist in
identifying the operations $m^{n}_{i}$ on $A$ and the operations $\mu^{n}_{i}$ on $Mag(A)$ which are deduced from the magmatic
operation $\cdot$\ .

The primitive part of a coassociative magmatic bialgebra is a
$\Prim Mag$-algebra by Proposition $\ref{MagV}$ and therefore it
is a $MagFine$ algebra by Proposition \ref{MuOperations}. Hence we
can apply the functor $U$ to $\Prim(\HH)$ for any $As^c\pt
Mag$-bialgebra $\HH$.

The second aim of this paper is to prove the structure theorem for
the triple of operads $(As, Mag, MagFine)$.

\begin{thm}\label{theorem2} Let $\HH$ be an $As^c\textrm{-}Mag$-bialgebra over the field $\KK$. Then, TFAE:\\
a) $\HH$ is connected (or conilpotent),\\ b) $\HH$ is isomorphic to $U(\Prim \HH)$,\\ c) $\HH$ is cofree in the category of connected
coassociative coalgebras.

In other words, $(As, Mag, MagFine)$ is a good triple of operads in the sense of \cite{Lod}.
\end{thm}

Observe that we did not make any characteristic assumption on $\KK$.

\begin{proo} The proof of this theorem is in two steps. First, we apply a result of \cite {Lod} to show that the structure theorem holds for the triple of operads $(As, Mag, \Prim Mag)$, cf.~Proposition \ref{structurePrim}. Second, we show that the morphism of operads $MagFine \to \Prim Mag$ is an isomorphism, cf.~Proposition \ref{isoFinePrim}.

In the Appendix  we give an alternative proof which is based on the construction of an explicit projection $\HH \epi \Prim \HH$.
\end{proo}

\begin{prop}\label{structurePrim} The structure theorem holds for the triple
$$(As, Mag, \Prim Mag).$$
\end{prop}
\begin{proo} In \cite {Lod} it is proved that for any $\CC^c\pt \AA$-bialgebra type whose compatibility is distributive and which satisfies the hypotheses (H1) and (H2) recalled below, the structure theorem holds.

The condition (H1) is the following:\\
the free $\AA$-algebra is equipped with a natural structure of $\CC^c\pt \AA$-bialgebra.

In Proposition \ref{MagV} we proved that $Mag(V)$ is a $As^c\pt Mag$-bialgebra.

The condition (H2) is the following:\\
the coalgebra map $\AA(V) \to \CC^c(V)$ has a natural splitting as coalgebra map.

Let us construct a coalgebra splitting of the map $Mag(V) \to As(V)$. We simply send the generator $\mu^c_n\in As_n$ to the left comb $\oo^n\in Mag_n$. It is immediate to check that it induces a coalgebra map $As^c(V) \to Mag(V)$.
\end{proo}
\begin{prop}\label{isoFinePrim} The morphism of operads $MagFine \to \Prim Mag$ constructed in section 2 is an isomorphism.
\end{prop}

\begin{proo} Since we showed in Proposition \ref{injective} that it is injective, it is sufficient to prove that the spaces $MagFine_n$ and $\Prim Mag_n$ have the same dimension. In \ref{Fineseries} we showed that $\dim MagFine_n = F_{n-1}$ (Fine number). Let us show that $\dim \Prim Mag_n = F_{n-1}$.

As a Corollary of Proposition \ref{structurePrim} we have an isomorphism of functors
$$Mag = As^c \circ \Prim Mag.$$
(Recall that the Poincar\'e series of a regular operad $\PP$ is $f_{\PP} (x) := \sum_{n\geq 1} \dim \PP_n x^n$).
From this isomorphism it follows that the Poincar\'e series of these functors are related by
$$f_{Mag}(x)= f_{As}(f_{\Prim Mag}(x)).$$
Since $f_{Mag}(x)= \frac{1-\sqrt{1-4x}}{2}$ and $ f_{As}(x) = \frac{x}{1-x}$, we get
 $$f_{\Prim Mag}(x) = \frac{1+2x-\sqrt{1-4x}}{2(2+x)}.$$
 Since this is precisely the generating series of the Fine numbers (up to a shift, cf. \cite {DS} and \ref{Fineseries}), we have $\dim \Prim Mag_n = F_{n-1}$.
\end{proo}

\begin{cor}
There is an equivalence of categories:
\begin{displaymath} \{\textrm{con.\ } As^c\pt Mag\bialg \}
{\begin{array}{c} U\\
\leftrightarrows\\
\Prim
\end{array}}
\{MagFine\alg \} \ .
\end{displaymath}
\end{cor}

\begin{cor} There is an isomorphism of functors
$$Mag = As \circ  MagFine\ .$$
\end{cor}

\subsection{Good triples of operads} In the terminology of \cite{Lod} the triple of operads $(As, Mag, MagFine)$ is \emph{good}, which means that the structure theorem holds. The interesting point is that it determines many other good triples, since, by a result of \cite{Lod}, it suffices to take some relation in $MagFine$ to construct a new good triple. For instance, if we mod out by the associator, then obviously the quotient of $MagFine$ becomes $Vect$ since all the higher operations are constructed from the associator. The middle operad becomes $As$ and we get the triple $(As, As, Vect)$ which has already been shown to be good in \cite{LR2}.

If we mod out by the pre-Lie relation
$$as(x,y,z)=as(x,z,y)$$
then the middle operad becomes the operad $preLie$ of pre-Lie algebras, and we get a good triple of operads
$$(As, preLie, MagFine / \cong).$$
From the isomorphism of functors $preLie = As^c \circ (MagFine / \cong )$ one can compute the dimension of $(MagFine / \cong )(n)$ from $n=1$ to $n=5$:
$$1, 0, 3, 16, 165.$$
It would be nice to find a small presentation of the
operad $MagFine / \cong$.

\bigskip
\bigskip

\section{Appendix}

{\bf Self-contained proof of the structure theorem for $(As, Mag, \Prim Mag)$ using an explicit idempotent.}

We first construct an explicit idempotent and we show some of its properties. Then we give a proof of the structure theorem.

\subsection{The idempotent}

Let $\HH$ be a connected $As^c\textrm{-}Mag$-bialgebra over the field $\KK$. We define an operator $e:{\overline {\HH}}\to {\overline {\HH}}$ as follows:

$$e = \id + \sum_{n\geq 1}(-1)^n \oleft{n+1}\circ \dd^n\ . $$

Note that the connectedness of $\HH$ implies that $e$ is well-defined.

For any family of elements $x_{1},\ldots ,x_{n}$ in $\HH$, let
$\oright{n}(x_{1},\ldots x_{n})$ be the right comb:
$$\oright{n}(x_{1},\ldots x_{n}):=x_{1}\cdot (x_{2}\cdot (\ldots
\cdot (x_{n-1}\cdot x_{n})\dots )).$$

\begin{lemma} The operator $e:{\overline {\HH}}\to {\overline {\HH}}$ has the following properties:\\
(i) $e(\oright{n}(x_1,\ldots,x_n))=0$ for all $x_1,\ldots,x_n\in\Prim\HH.$\\
 (ii) $\dd\circ e = 0.$\\
  (iii) $e\vert_{\Prim\HH} = \id_{\Prim\HH}.\\$
\end{lemma}

\begin{proo} (i): Note that if the elements $x_{1},\dots ,x_{n}$ belong to $Prim(\HH)$, then:
$$\dd ^{k}(\oright{n}(x_1,\ldots,x_n))=\sum _{1\leq i_{1}< \dots <i_{k}\leq n}\oright{i_{1}}(x_1,\ldots,x_{i_{1}})\otimes \dots \otimes 
\oright{n-i_{k}}(x_{i_{k}+1},\ldots,x_n) .$$
So, $$\displaylines{
\oleft{k+1}\circ \dd ^{k}(\oright{n}(x_1,\ldots,x_n))=\hfill \cr
\sum \oleft{k+1}(\oright{i_{1}}(x_1,\ldots,x_{i_{1}}), \dots ,
 \oright{n-i_{k}}(x_{i_{k}+1},\ldots,x_n))=\cr
 \sum \oleft{k+2}(x_{1},\oright{i_{1}-1}(x_2,\ldots,x_{i_{1}}),\dots ,
 \oright{n-i_{k}}(x_{i_{k}+1},\ldots,x_n))=\cr
  \sum \oleft{k+2}(\oright{1}(x_{1}),\oright{i_{1}-1}(x_2,\ldots,x_{i_{1}}),\dots , \oright{n-i_{k}}(x_{i_{k}+1},\ldots,x_n)).\cr}$$
Each term $\oleft{k+1}(\oright{i_{1}}(x_1,\ldots,x_{i_{1}}), \dots ,\oright{n-i_{k}}(x_{i_{k}+1},\ldots,x_n))$ appears twice in $e(\oright{n}(x_1,\ldots,x_n))$:
in $\oleft{k+1}\circ \dd ^{k}(\oright{n}(x_1,\ldots,x_n))$ with coefficient $(-1)^{k}$, and in 
$\oleft{k+2}\circ \dd ^{k+1}(\oright{n}(x_1,\ldots,x_n))$ with coefficient $(-1)^{k+1}$. Since $e(\oright{n}(x_1,\ldots,x_n))$ is the sum of these terms, we get that 
$e(\oright{n}(x_1,\ldots,x_n))=0$.

(ii): Using the unital infinitesimal relation verified by $\dd$ it is immediate to check that, 
for elements $x_{1},\dots ,x_{n}$ of ${\overline {\HH}}$, we have the following equality:
$$\displaylines {
\dd \circ \oleft{n}(x_{1},\dots ,x_{n})=\sum _{i=1}^{n-1}\oleft{i}(x_{1},\dots ,x_{i})\otimes \oleft{n-i}(x_{i+1},\dots ,x_{n})\cr
\hfill + \sum _{j=1}^{n}\oleft{j}(x_{1},\dots ,x_{j-1},x_{j(1)})\otimes \oleft{n-j+1}(x_{j(2)},x_{j+1},\dots ,x_{n}),\cr }$$
where $\delta (x_{j})=x_{j(1)}\otimes x_{j(2)}$.
The coassociativity of $\dd$ and the formula above imply that:
$$\displaylines{\dd \circ \oleft{n}\circ \dd ^{n-1}(x)=\sum _{i=1}^{n-1}\oleft{i}(x_{(1)}, \dots ,x_{(i)})\otimes \oleft{n-i}(x_{(i+1)}, \dots ,x_{(n)})+\cr
\hfill \sum _{i=1}^{n}\oleft{i}(x_{(1)},\dots ,x_{(i)})\otimes \oleft{n+1-i}(x_{(i+1)},\dots ,x_{(n+1)}).\cr }$$
So, for $i\geq 1$, the term $\oleft{i}(x_{(1)}, \dots ,x_{(i)})\otimes \oleft{n-i}(x_{(i+1)}, \dots ,x_{(n)})$ appears twice in $\dd\circ e (x)$. 
In $\dd \circ \oleft{n}\circ \dd ^{n-1}(x)$ with coefficient $(-1)^{n-1}$, and in $\dd \circ \oleft{n-1}\circ \dd ^{n-2}(x)$ 
with coefficient $(-1)^{n-2}$, which ends the proof.

(iii): The assertion follows directly from the definition of $e$ and the definition of $\dd^n$.

\end{proo}
\medskip

\begin{lemma}\label{lemafinal} Given a connected $As^c\textrm{-}Mag$-bialgebra $\HH$, for any $x\in {\overline {\HH}}$ one has the following equality:
$$x=e(x)+\oright{2}((e\t e)( \dd ^{1}(x)))+\ldots
+\oright{n+1}(e^{\t n+1}(\delta ^{n}(x)))+\ldots .$$
\end{lemma}
\medskip

\begin{proo} Let $x\in F_{r}\HH$, with $r\geq 1$. If $r=1$, then the result is obvious.

For $r>1$, note that $e(x)=x - e(x_{(1)})\cdot x_{(2)}$, and
that $$\oright{n+1}(e^{\t n+1}(\delta ^{n}(x)))=
e(x_{(1)})\cdot \oright{n}(e^{\t n}(\delta ^{n-1}(x_{(2)})),$$
where $\dd (x)=  x_{(1)}\t x_{(2)}$ (here we omit the summation sign). Since any $x_{(2)}\in
F_{r-1}\HH$, it verifies the equality, so one gets:
\begin{eqnarray*} \sum _{i\geq 0}\oright{i+1}(e^{\t i+1}(\dd
^{i}(x)))&=& x-e(x_{(1)})\cdot x_{(2)}+\sum
_{i\geq 1}\oright{i+1}(e^{\t i+1}(\dd ^{i}(x)))\\
&=& x-
e(x_{(1)})\cdot x_{(2)}+\sum _{i\geq 1}(e(x_{(1)})\cdot
\oright{i}(e^{\t i}(\delta ^{i-1}(x_{(2)})))\\
&=& x-
e(x_{(1)})\cdot x_{(2)}+ e(x_{(1)})\cdot x_{(2)}=x.\\
\end{eqnarray*}

\end{proo}

Let us now prove  Theorem \ref{theorem2}.

\begin{proo} $a)\Rightarrow b)$
For a connected coassociative magmatic bialgebra $\HH$ we consider the map
$\varphi :\HH\to U(\Prim \HH)$ given by the composition of the
projection $p: Mag(\Prim \HH)\rightarrow U(\Prim \HH)$ with the
map $$x\mapsto \sum _{n\geq 1}\big( (\oright{n};e(x_{(1)})\ldots
e(x_{(n)}))\big),$$ where $\dd ^{n-1}(x)= x_{(1)}\t \ldots \t
x_{(n)}$.

It is not difficult to check that $\varphi$ is a coalgebra morphism.

To prove that it is a magmatic algebra morphism, we first verify the formula $\varphi(x\cdot y)= \varphi(x)\cdot \varphi(y)$ for $x \in \Prim \HH$ and any $y$. Second, we prove, by a recursive argument on $n$, that this formula holds also for 
$x= (\oright{n} ; x_1, \ldots, x_n)$ and $y= (\oright{m} ; y_1, \ldots, y_m)$, $x_i, y_j \in \Prim \HH$.  Lemma \ref{lemafinal} implies that the  formula holds also for all $x$ and $y$.

To construct the inverse map  we note that there
exists a unique homomorphism of magmatic algebras $\epsilon:Mag(\Prim
\HH)\rightarrow\HH$, which extends the canonical injection of $\iota :\Prim \HH \hookrightarrow \HH$. For
$n\geq 3$, given elements $x_{1},\ldots ,x_{n}$ in $\Prim \HH$ and
an integer $1\leq i\leq n-2$, let $\mu _{i}^{n}(x_{1},\ldots
,x_{n})$ be the image of $x_{1}\t \ldots \t x_{n}$ under the
operation $\mu _{i}^{n}$ of $\Prim \HH$ and let ${\overline {\mu
}}_{i}^{n}(x_{1},\ldots ,x_{n})$ be the image of $x_{1}\t \ldots
\t x_{n}$ under the operation ${\overline {\mu}} _{i}^{n}$ which
is constructed using the magmatic product of $Mag(\Prim \HH)$.

\noindent Since $$\epsilon(\mu _{i}^{n}(x_{1},\ldots
,x_{n}))= {\mu }_{i}^{n}(x_{1},\ldots ,x_{n})=
\epsilon ({\overline {\mu }}_{i}^{n}(x_{1},\ldots ,x_{n})),$$ for any
$x_{1},\ldots ,x_{n}\in \Prim \HH$, there exists a unique
homomorphism $\psi :U(\Prim \HH)\rightarrow \HH$, which factors $\epsilon $.

Lemma \ref{lemafinal} implies that $\psi \circ \varphi
=id_{\HH}$. To see that $\varphi \circ \psi =id_{U(\Prim
\HH)}$, it suffices to note that $\varphi \circ \epsilon $ and  $p:
Mag(\Prim \HH)\rightarrow U(\Prim \HH)$  extend the composition $\Prim \HH\hookrightarrow \HH \xrightarrow
{\varphi } U(\Prim \HH)$ to $Mag(\Prim \HH)$, so $\varphi \circ \epsilon =p$. This equality implies that $id_{U(\Prim \HH)}=\varphi
\circ \psi $.
\medskip

$a)\Rightarrow c)$ We note that the map $x\mapsto \alpha (x):=\sum
_{n\geq 1}(\sum e(x_{(1)})\t \ldots \t e(x_{(n)}))$ gives a linear
bijection from $\HH $ into $T^{c}(\Prim \HH )$. Moreover, 
$$\Delta^c(e(x_{(1)})\t \ldots \t e(x_{(n)}))=\sum
_{i=0}^{n}(e(x_{(1)})\t \ldots \t e(x_{(i)}))\t (e(x_{(i+1)}) \t
\ldots \t e(x_{(n)})),$$ for $n\geq 1$, where $\Delta ^{c}$ is the
deconcatenation coproduct on $T^{c}(\Prim \HH )$. On the other
hand, let $\delta ^{i-1}(x_{(1)})=\sum x_{(11)}\t \ldots \t
x_{(1i)}$ and $\delta ^{n-i-1} (x_{(2)})= \sum x_{(2,i+1)}\t
\ldots \t x_{(2n)}$, for $\Delta (x)=\sum x_{(1)}\t x_{(2)}.$ The
coassociativity of $\Delta $ implies that 
$$\delta ^{n-1}(x)=\sum
x_{(1)}\t \ldots \t x_{(n)}=\sum (\sum x_{(11)}\t \ldots \t
x_{(2n)}).$$
Then we have
$$\displaylines { \sum \alpha (x_{(1)})\t \alpha
(x_{(2)})=\hfill \cr \sum (\sum _{i\geq 1}e(x_{(11)})\t \ldots \t
e(x_{(1i)}))\t (\sum _{n-i\geq 1}e(x_{(2,i+1)}\t \ldots \t
e(x_{(2n)})))=\cr \sum _{n\geq 1}(\sum
_{i=0}^{n}((e(x_{(1)})\t\ldots \t e(x_{(i)}))\t
(e(x_{(i+1)})\t\ldots \t e(x_{(n)})) )=\Delta
^{c}(\alpha (x)),\cr }$$ which implies that $\alpha $ is a
coalgebra isomorphism.
\medskip

$c)\Rightarrow b)$ We consider the coalgebra isomorphism $
T^{c}(\Prim \HH )\rightarrow U(\Prim \HH)$ which sends $x_{1}\t
\ldots \t x_{n} \mapsto \oright{n}(x_{1},\ldots ,x_{n})$. Since
$\HH $ is cofree, we obtain that $\HH$ is isomorphic to $U(\Prim
\HH)$ as a coalgebra. An easy calculation shows that the
isomorphism is also an isomorphism of magmatic algebras.

$b)\Rightarrow a)$ is obvious.
\end{proo}

\bigskip

 \noindent {\bf Acknowledgement.}

 R.~Holtkamp wants to thank for the hospitality of the I.R.M.A.(Universit\'e Louis Pasteur, Strasbourg). M.~Ronco acknowledges the financial aid of the program FONDECYT 1060224.
 R.~Holtkamp and J.-L.~Loday thank the Mittag-Leffler Institute for an invitation in April 2004 during which this work was begun. The authors thank Bruno Vallette for his remarks on a first version of this paper.


\end{document}